\documentclass[10pt,twoside,twocolumn]{IEEEtran}
\usepackage[latin9]{inputenc}
\usepackage{float}
\usepackage{amsmath}
\usepackage{amsthm}
\usepackage{amssymb}
\usepackage{graphicx}

\makeatletter

\providecommand{\tabularnewline}{\\}
\floatstyle{ruled}
\newfloat{algorithm}{tbp}{loa}
\providecommand{\algorithmname}{Algorithm}
\floatname{algorithm}{\protect\algorithmname}

\theoremstyle{plain}
\newtheorem{thm}{\protect\theoremname}
\theoremstyle{plain}
\newtheorem{lem}[thm]{\protect\lemmaname}
\theoremstyle{remark}
\newtheorem{rem}[thm]{\protect\remarkname}
\theoremstyle{definition}
\newtheorem{defn}[thm]{\protect\definitionname}
\theoremstyle{plain}
\newtheorem{cor}[thm]{\protect\corollaryname}

\usepackage{array}\usepackage{verbatim}\usepackage{amsfonts}\usepackage{color}

\usepackage{enumerate}
\usepackage{epstopdf}
\usepackage[noadjust]{cite}
  
\usepackage{dblfloatfix}
\usepackage{dsfont}
\usepackage{breqn}

\usepackage{url}

\newcommand{\R}{\mathbb{R}}			
\newcommand{\C}{\mathbb{C}}			

\newcommand{\I}{\mathcal{I}}

\newcommand{\tr}{\tilde{\rho}_2}
\DeclareMathOperator{\ee}{\mathbb{E}}			

\newcommand{\re}[1]{\text{Re}\left(#1\right)} 
\newcommand{\img}[1]{\text{Im}\left(#1\right)} 

\floatstyle{ruled}
\newfloat{algorithm}{tbp}{loa}
\providecommand{\algorithmname}{Algorithm}
\floatname{algorithm}{\protect\algorithmname}

\usepackage{amsfonts}

\providecommand{\corollaryname}{Corollary}
\providecommand{\definitionname}{Definition}
\providecommand{\theoremname}{Theorem}

\global\long\def\ee{\mathbb{E}}
\global\long\def\C{\mathbb{C}}
\global\long\def\I{\mathcal{I}}
\global\long\def\tr{\mathrm{tr}}
\global\long\def\K{\mathcal{K}}
\global\long\def\one{\mathbf{1}}
\global\long\def\vect{\mathrm{vec}}

\usepackage{tikz}
\usepackage{circuitikz}
\usepackage{tikz-cd}
\usetikzlibrary{shapes}
\tikzstyle{block} = [draw,rectangle, rounded corners, minimum width=1cm, minimum height=0.8cm,text centered, line width=1pt ]
\tikzstyle{agent} = [draw,circle, minimum width=1.1cm, text centered, line width=1pt ]
\tikzstyle{arrow} = [thick,->,>=stealth,line width=1pt]
\usetikzlibrary{arrows,decorations.markings,decorations}
\tikzset{
    addarrow/.style={decoration={markings, mark=at position 1 with {\arrow{stealth}}},
                     postaction={decorate}}
}

\@ifundefined{showcaptionsetup}{}{%
 \PassOptionsToPackage{caption=false}{subfig}}
\usepackage{subfig}
\makeatother

\providecommand{\corollaryname}{Corollary}
\providecommand{\definitionname}{Definition}
\providecommand{\lemmaname}{Lemma}
\providecommand{\remarkname}{Remark}
\providecommand{\theoremname}{Theorem}

\begin{document}
\title{Large-Scale Traffic Signal Offset Optimization }
\author{Yi Ouyang\thanks{Y. Ouyang is with Preferred Networks, Inc. (email: ouyangyi@preferred-america.com)},
Richard Y. Zhang\thanks{R. Y. Zhang is with the Department of Electrical and Computer Engineering,
University of Illinois at Urbana\textendash Champaign (email: ryz@illinois.edu).}, Javad Lavaei\thanks{J. Lavaei is with the Department of Industrial Engineering and Operations
Research and also the Tsinghua-Berkeley Shenzhen Institute, University
of California, Berkeley (email: lavaei@berkeley.edu).}, and Pravin Varaiya\thanks{P. Varaiya is with the Department of Electrical Engineering and Computer
Sciences, University of California, Berkeley (varaiya@berkeley.edu).} \thanks{This work was supported by the ONR grant N00014-17-1-2933, DARPA grant
D16AP00002, AFOSR grant FA9550- 17-1-0163, ARO grant W911NF-17-1-0555,
and NSF Awards 1545116 and 1807260. A preliminary and abridged version
has appeared in the Proceedings of the 57th IEEE Conference on Decision
and Control, 2018 \cite{ouyang2018cdc}. } }

\maketitle

\begin{abstract}
The offset optimization problem seeks to coordinate and synchronize
the timing of traffic signals throughout a network in order to enhance
traffic flow and reduce stops and delays. Recently, offset optimization
was formulated into a continuous optimization problem without integer
variables by modeling traffic flow as sinusoidal. In this paper, we
present a novel algorithm to solve this new formulation to near-global
optimality on a large-scale. Specifically, we solve a convex relaxation
of the nonconvex problem using a tree decomposition reduction, and
use randomized rounding to recover a near-global solution. We prove
that the algorithm always delivers solutions of expected value at
least 0.785 times the globally optimal value. Moreover, assuming that
the topology of the traffic network is ``tree-like'', we prove that
the algorithm has near-linear time complexity with respect to the
number of intersections. These theoretical guarantees are experimentally
validated on the Berkeley, Manhattan, and Los Angeles traffic networks.
In our numerical results, the empirical time complexity of the algorithm
is linear, and the solutions have objectives within 0.99 times the
globally optimal value. 
\end{abstract}


\IEEEpeerreviewmaketitle{}

\section{Introduction}

In transportation engineering, \emph{traffic signal timing} is the
problem of selecting and adjusting the timing of traffic lights in
order to reduce congestion and improve traffic flow. This classical
problem is commonly formulated as three subproblems:
\begin{itemize}
\item \emph{Cycle length optimization}, where the total network is divided
into subsections, and a common cycle period is assigned to each subsection;
\item \emph{Green split optimization}, where traffic lights within the same
intersection are timed to improve throughput; and
\item \emph{Offset optimization}, where traffic lights over different intersections
are coordinated to enhance network-wide performance.
\end{itemize}
Ideally, these subproblems would be solved simultaneously for the
best performance~\cite{gartner1975optimization2,gartner1976mitrop}.
Owing to issues of computational tractability, however, the established
practice is an \emph{iterative} procedure: manually divide the network
into subsections, sweep the cycle length over a range of values, and
solve the green split and offset optimization subproblems alternatingly
for each fixed cycle length~\cite{gartner1973generalized,gartner1975optimization}.
This is precisely the solution procedure implemented in the industry-standard
software packages TRANSYT-7F~\cite[Sec 2.4]{wallace1984transyt}
and Synchro~\cite[Ch. 18]{synchro2017}.

In this paper, we focus our attention on the offset optimization subproblem.
The goal is to create \emph{green waves}, in which green lights are
synchronized to allow a car to drive through multiple intersections
without stopping for a red light, and to maximize the length or \emph{bandwidth}
of these green waves. Clearly, green waves are only possible if cycle
lengths are the same, or else the synchronization would be lost over
time. For this reason, the standard model represents traffic flow
as square waves with a common cycle length but separate green times
and red times. The exact green splits are assumed to be given and
fixed, with the understanding that they will be separately optimized
at a later stage. 

\subsection{Previous Approaches}

The offset optimization problem is highly nonconvex, so solution approaches
based on incremental adjustments\textemdash such as those implemented
in TRANSYT and Sychro\textemdash can get stuck at a locally optimal
solution. In order to obtain a \emph{globally} optimal solution, the
standard approach is to reformulate the problem into a mixed-integer
program~\cite{little1966synchronization,little1981maxband,gartner2002arterial,gartner2004progression}
and apply a general-purpose integer programming solver like Gurobi
or CPLEX. The latter approach is highly effective on a small scale,
but\textemdash as is typical for techniques based around integer programming\textemdash suffers
from severe computational issues as the problem size grow large. In
practice, it may not even find a feasible point that does not violate
constraints in a reasonable amount time, let alone a globally optimal
solution.

Instead, computing globally optimal solutions to large-scale networks
generally requires simplifying assumptions. In particular, if a penalty
function known as a \emph{link delay function} is assigned to each
road link with respect to the offset difference, then dynamic programming
can be used to minimize the sum of all link delay functions~\cite{allsop1968selection,gartner1972optimal,gartner1973generalized}.
For certain network topologies, this approach is guaranteed to compute
a globally optimal solution in linear time. However, it is often tricky
to choose a link delay function that accurately reflects real-world
considerations like queues, delays, and green waves~\cite{allsop1968selection,gartner1975optimization}.
Also, its use relies on an assumption of link independence that may
not be fully realistic~\cite{gartner1975optimization}.

Recently, Coogan et al.~\cite{coogan2015offset,coogan2017offset}
proposed an approach, based on an alternative simplifying assumption,
that outperforms the link delay function approach described above~\cite{coogan2017offset}
and the incremental adjustment approach found in Synchro~\cite{amini2018optimizing}.
By modeling traffic flow as \emph{sinusoidal}, the problem of minimizing
total queue lengths can be posed as a quadratically-constrained quadratic
program (QCQP). The QCQP is nonconvex, but can be relaxed into a convex
semidefinite program (SDP) using standard techniques, and solved using
an interior-point method. In turn, the solution to the SDP can often
recover a globally optimal solution for the QCQP. If desired, the
solution can be further refined using TRANSYT or Synchro~\cite{kim2017offset}.

Nevertheless, the Coogan et al.~\cite{coogan2015offset,coogan2017offset}
approach suffers two serious computational issues that prevent its
use on real traffic networks. First, the approach often yields, but
does not guarantee, a globally optimal solution. Indeed, such a guarantee
is not even possible in general unless P=NP. Moreover, the convex
SDP that underpins the approach has a worst-case solution complexity
of $O(n^{4.5})$ time and $O(n^{2})$ memory. While these figures
are formally polynomial, their large exponents limit the number of
intersections $n$ to no more than a few hundred.

\subsection{Main Results}

Our main contribution in this paper is an algorithm that is \emph{guaranteed}
to solve the formulation of Coogan et al.~\cite{coogan2017offset}
to \emph{near-global} optimality in \emph{near-linear} time. In Section~\ref{sec:Approximation-Algorithm},
we prove that the algorithm always delivers solutions of expected
value at least $\pi/4\ge0.785$ times the globally optimal value.
Moreover, assuming that the topology of the traffic network is ``tree-like'',
we prove in Section~\ref{sec:Efficient-Implementation-for} that
the algorithm has near-linear $O(n^{1.5})$ time complexity and linear
$O(n)$ memory complexity with respect to the number of intersections
$n$. These theoretical guarantees are experimentally validated in
Section~\ref{sec:numerical} on the Berkeley, Manhattan, and Los
Angeles traffic networks. In our numerical results, the algorithm
achieves a linear empirical time complexity, and the solutions found
all have objectives within 0.99 times the globally optimal value. 

Our algorithm works by reformulating offset optimization into a complex-valued
quadratically-constrained quadratic program (QCQP) with a similar
form to the classic MAX-CUT problem in combinatorial optimization
\cite{karp1972reducibility}, and relaxing the QCQP into a semidefinite
program (SDP). Inspired by the Goemans\textendash Williamson algorithm
for MAX-CUT \cite{goemans1995improved}, we prove that projecting
the SDP solution onto a random hyperplane recovers a solution to the
QCQP with an approximation ratio of $\pi/4$. We solve the SDP relaxation
using the sparsity-exploiting chordal conversion technique of Fukuda
et al.~\cite{fukuda2001exploiting} and the dualization technique
recently developed by Zhang and Lavaei~\cite{zhang2017sparse,zhang2017sparse2}.
Directly solving the SDP in the complex domain yields significant
improvement on runtime, compared to our previous results in the conference
version of this paper~\cite{ouyang2018cdc}. When a network is ``sparse''
in the sense that it has a bounded treewidth \cite{robertson1986graph},
we prove that the overall algorithm has worst-case complexity of $O(n^{1.5})$
time and $O(n)$ memory. 

\subsection{Future work}

Our algorithm solves the Coogan et al.~\cite{coogan2017offset} formulation
of offset optimization to near-global optimality and in near-linear
time. The corresponding traffic model assumes sinusoidal traffic flow
that share a common global cycle length, but this may not always be
realistic under all traffic conditions. There have been efforts to
validate the realism of the model, but full-scale studies have yet
to be performed, owing largely to the lack of efficient algorithms.
Our algorithm makes it possible to solve the formulation on a large
scale, so an important next step is to validate the traffic model,
by benchmarking our results against realistic large-scale micro-simulations.

Also, our algorithm optimizes the offsets on their own, with the understanding
that the cycle lengths and green splits would be separately optimized
at a later time, possibly in an alternating fashion with respect to
the offsets. Nevertheless, a truly globally optimal signal timing
profile would require all three parameters to be simultaneously coordinated.
It remains future work to benchmark the performance of the alternating
global optimization approach against a simultaneous local optimization
approach.

\subsection*{Notation}

The sets $\R$ and $\C$ are the real and complex numbers. Subscripts
indicate element-wise indexing. The notation $X_{\mathcal{I},\mathcal{J}}$
indicates the submatrix of $X$ indexed by columns sets $\mathcal{I},\mathcal{J}\subseteq\{1,2,\dots,n\}$.
The superscripts ``$T$'' and ``$H$'' refer to the transpose
and the Hermitian transpose. We write $i=\sqrt{-1}$ and use $\text{Re}(x)$,
$\text{Im}(x)$, $\bar{x}$, $\angle x$, and $|x|$ to denote the
real part, imaginary part, conjugate, angle, and absolute value. The
identity matrix is $I$ and the vector-of-ones is $\one$; their sizes
are inferred from context. The trace, rank, and column vectorization
are denoted $\tr(X)$, $\text{rank}(X)$, and $\vect(X)$. $X\succeq0$
means that $X$ is Hermitian and positive semidefinite. $|\mathcal{S}|$
denotes the cardinality of a set $\mathcal{S}$.

\section{Problem Formulation}

\label{sec:formulation}

To determine traffic signal offsets, we adopt the traffic network
model with sinusoidal approximation proposed in \cite{coogan2017offset}.
In what follows, we will first describe the model and explain this
sinusoidal approximation technique. Then, using this model, we formulate
a mathematical optimization problem to select offsets that minimize
the lengths of vehicle queues of the networks.

\subsection{Traffic Network Model}

\begin{figure}[H]
\centering \begin{tikzpicture}[scale=0.9]
\node[agent] (1) at  (0,0) {\includegraphics[width = 1.1cm]{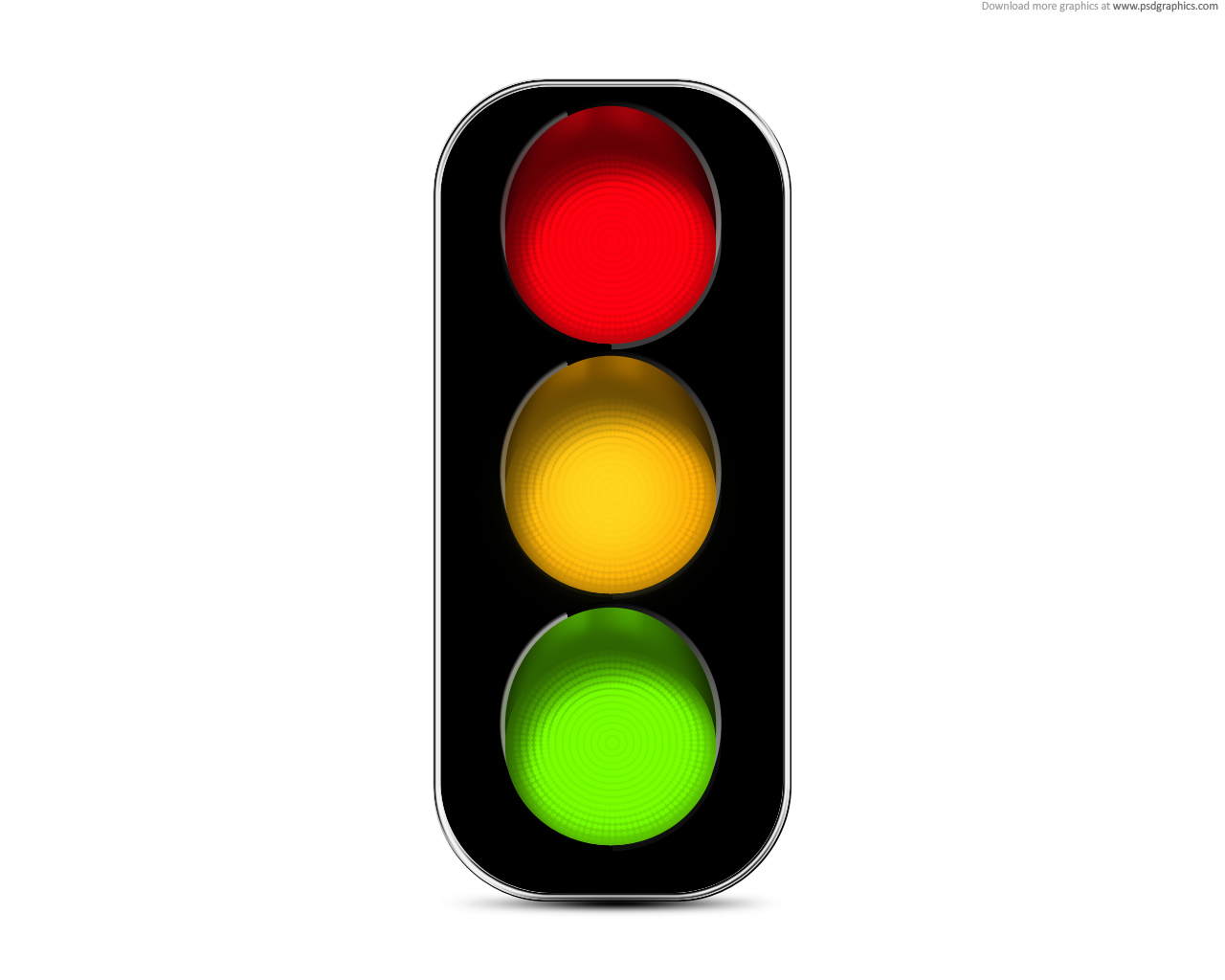}};
\node at (-0.5,0) {$\tau(l)$};
\node[agent] (2) at  (3,0) {\includegraphics[width = 1.1cm]{traffic-signal}};
\node[agent] (3) at  (0,-3) {\includegraphics[width = 1.1cm]{traffic-signal} };
\node at (-0.5,-3) {$\sigma(l)$};
\node[agent] (4) at  (3,-3) {\includegraphics[width = 1.1cm]{traffic-signal}};
\node[agent] (5) at  (-3, -3) {Source $\epsilon$};
\draw[arrow] (1.10) to (2.170);
\draw[arrow] (2.190) to (1.350);
\draw[arrow] (1.260) to node[left]{$l$} (3.100);
\draw[arrow] (1.325) to (4.125);
\draw[arrow] (2) to (4);
\draw[arrow] (3.80) to (1.280);
\draw[arrow] (3.10) to (4.170);
\draw[arrow] (4.190) to (3.350);
\draw[arrow] (4.145) to (1.305);
\draw[arrow] (5) to (3);
\draw[arrow] (5) to (1);
\end{tikzpicture} \caption{Traffic Network}
\end{figure}

Consider a traffic network described by a directed graph $G=(\mathcal{S}\cup\{\epsilon\},\mathcal{L})$.
Each node of the graph represents an intersection; node $i\in\mathcal{S}=\{1,2,\dots,|S|\}$
represents a signalized intersection and node $\epsilon$ is the dummy
intersection (source) for traffic originating outside the network.
Let $n=|S|+1$ be the number of intersections including the dummy
intersection. The dummy node $\epsilon$ is also referred to as node
$n$. Each directed edge in $\mathcal{L}$ represents a traffic link
between two intersections/signals and the vehicle queue associated
with the link. For each $l\in\mathcal{L}$, $\tau(l)\in\mathcal{S}$
indicates its upstream intersection and $\sigma(l)\in\mathcal{S}$
represents the downstream intersection which serves the queue of the
link. $\mathcal{E}=\{l\in\mathcal{L},\tau(l)=\epsilon\}\subset\mathcal{L}$
is the set of entry links that direct exogenous traffic from the dummy
intersection (source) to the network; other links are non-entry links
and the travel time from its upstream to downstream intersections
is denoted by $\lambda_{l}$. There is no need to explicitly model
links that exit the network because exiting traffic are considered
in the calculation of turn ratios, which will be defined later.

The vehicle queue associated with each link $l\in\mathcal{L}$ has
length $q_{l}(t)$ at time $t$. The queue length $q_{l}(t)$ follows
a continuous-time fluid queue model given by 
\begin{align}
\dot{q}_{l}(t)=a_{l}(t)-d_{l}(t)\label{eq:q_def}
\end{align}
where $a_{l}(t)$ is the arrival rate for vehicles arriving from the
upstream intersection and $d_{l}(t)$ is the departure rate that depends
on the downstream intersection signal. Both $a_{l}(t)$ and $d_{l}(t)$
are in units of vehicles per hour.

Vehicles coming from a link are allowed to pass through an intersection
when the link is activated by the traffic signal, i.e., green light
for the link. To avoid collision, each signal switches among activation
patterns of non-conflicting links according to a signal control sequence.
All intersections are assumed to operated under fixed time control
\cite{muralidharan2015analysis} with common cycle. This means that
the signal control sequence of each intersection has a fixed periodic
cycle, and all intersections have a common cycle time $T=1$ time
unit.

The signal offset $\theta_{s}\in[0,1)$ for an intersection $s\in\mathcal{S}$
represents the phase difference of the signal control sequence from
a global clock. 
For each link $l\in\mathcal{L}$, vehicles from its queue is allowed
to pass through intersection $\sigma(l)$ at times $n+\theta_{\sigma(l)}+\gamma_{l}$
for $n=0,1,2,\ldots$, where $\gamma_{l}\in[0,1)$ is called the link's
green split that represents the time difference of the midpoint of
the activation time for the link and the beginning of the offset time
$\theta_{\sigma(l)}$. For $l,k\in\mathcal{L}$, the turn ratio $\beta_{lk}\in[0,1]$
denotes the fraction of vehicles that are routed to link $k$ upon
exiting link $l$. When $\sigma(l)\neq\tau(k)$, $\beta_{lk}=0$ because
the two links are not connected. For every link $l\in\mathcal{L}$
it holds that 
\begin{align*}
\sum_{k\in\mathcal{L}}\beta_{lk}\leq1
\end{align*}
where strict inequality in the above equation models the situation
that a fraction of vehicles exit the network via an unmodeled link
from intersection $\sigma(l)$.

Similarly to \cite{coogan2017offset}, we assume that the network
is in the periodic steady state and approximate all arrivals, departures,
and queue lengths by sinusoid functions with period $T=1$. Specifically,
the departure rate of link $l$ is assumed to be 
\begin{align*}
d_{l}(t)= & f_{l}(1+\cos(2\pi(t-\theta_{\sigma(l)}-\gamma_{l})))
\end{align*}
where $f_{l}$ is the average departure rate of link $l$. By defining
$z_{j}=e^{i2\pi\theta_{j}}$ for $j\in\mathcal{S}$ and $D_{l}=f_{l}e^{-i2\pi\gamma_{l}}$,
one can write the departure rate at link $l$ as 
\begin{align}
d_{l}(t)= & f_{l}+\re{e^{i2\pi t}D_{l}\bar{z}_{\sigma(l)}}.\label{eq:departure_rate}
\end{align}
Since vehicles arrive at a non-entry link from its upstream links
after a delay equal to the travel time, the arrival rate of a non-entry
link $l\in\mathcal{L}\setminus\mathcal{E}$ is given by 
\begin{align*}
a_{l}(t)=\sum_{k\in\mathcal{L}}\beta_{kl}d_{k}(t-\lambda_{l}).
\end{align*}
The periodic steady-state assumption implies that the average arrival
rate is the same as the average departure rate at each link \cite{muralidharan2015analysis},
i.e., 
\begin{align*}
\int_{0}^{1}a_{l}(t)dt=\int_{0}^{1}d_{l}(t)dt.
\end{align*}
Therefore, we have 
\begin{align*}
\sum_{k\in\mathcal{L}}\beta_{kl}f_{k}=f_{l}.
\end{align*}
Then, the arrival rate can be further expressed as 
\begin{align}
a_{l}(t)= & f_{l}+\re{e^{i2\pi t}A_{l}\bar{z}_{\tau(l)}}
\label{eq:arrival_rate}
\end{align}
where $A_{l}=e^{-i2\pi\lambda_{l}}\sum_{k\in\mathcal{L}}\beta_{kl}D_{k}$.

For an entry link $l\in\mathcal{E}$, the approximation assumes that
\begin{align}
a_{l}(t) & =f_{l}+\alpha_{l}\cos(2\pi(t-\phi_{l})))\nonumber \\
 & =f_{l}+\re{e^{i2\pi t}A_{l}\bar{z}_{\tau(l)}}\label{eq:arrival_rate_entry}
\end{align}
where $z_{\tau(l)}=e^{i2\pi\theta_{n}}=1$ with the offset $\theta_{n}$
of the dummy intersection $\epsilon$ (intersection $n$) defined
to be $0$ in the above equation, $\alpha_{l}\leq f_{l}$ is the relative
amplitude of the arrival peak minus the average rate, $A_{l}=\alpha_{l}e^{-2\pi\phi_{l}}$,
and $\phi_{l}\in[0,1)$ is the offset for the center of the arrival
peak.

It follows from the queue dynamics \eqref{eq:q_def}, departure rate
\eqref{eq:departure_rate} and arrival rate \eqref{eq:arrival_rate}-\eqref{eq:arrival_rate_entry}
of the links that the queue length $q_{l}(t)$ of each link $l\in\mathcal{L}$
evolves according to the equation 
\begin{align*}
\dot{q}_{l}(t) & =a_{l}(t)-d_{l}(t)\\
 & =\re{e^{i2\pi t}(A_{l}\bar{z}_{\tau(l)}-D_{l}\bar{z}_{\sigma(l)})}.
\end{align*}
Accordingly, the average queue length at link $l$, denoted by $Q_{l}$,
is given by 
\begin{align*}
Q_{l}=\frac{1}{2\pi}|(A_{l}\bar{z}_{\tau(l)}-D_{l}\bar{z}_{\sigma(l)})|.
\end{align*}

\subsection{Offset Optimization Problem}

\label{subsec:formulation} The average queue lengths $Q_{l}$ where
$l\in\mathcal{L}$, are important performance metrics for traffic
networks. Following the approach in \cite{coogan2017offset}, we formulate
the offset optimization problem as selecting offsets $\theta_{s},s=1,2,\dots,n$
with the goal of minimizing the total average squared queue length.
Note that the queue lengths are invariant to a constant shift for
all $\theta_{s}$ where $s=1,2,\dots,n$. Therefore, instead of restricting
$\theta_{n}=0$ for the dummy intersection $\epsilon$, one can allow
$\theta_{n}$ to be a variable that takes any value in the interval
$[0,1)$ and set the offset of each intersection $s\in\mathcal{S}$
to be the relative offset $\theta_{s}-\theta_{n}$. Then, the offset
optimization problem can be formulated as follows: 
\begin{align}
\underset{\theta_{1},\dots,\theta_{n}}{\text{minimize }} & \sum_{l\in\mathcal{L}}Q_{l}^{2}\label{eq:problem}\\
\text{subject to } & Q_{l}=\frac{1}{2\pi}|(A_{l}\bar{z}_{\tau(l)}-D_{l}\bar{z}_{\sigma(l)})|\nonumber \\
 & z_{s}=e^{i2\pi\theta_{s}},\quad s=1,2,\dots,n.\nonumber 
\end{align}
Note that the queue length of each link satisfies 
\begin{align*}
Q_{l}^{2}= & \frac{1}{(2\pi)^{2}}|(A_{l}\bar{z}_{\tau(l)}-D_{l}\bar{z}_{\sigma(l)})|^{2}\\
= & \frac{1}{(2\pi)^{2}}(|A_{l}|+|D_{l}|)^{2}\\
- & \frac{1}{(2\pi)^{2}}(2|A_{l}||D_{l}|+\bar{D_{l}}A_{l}\bar{z}_{\tau(l)}z_{\sigma(l)}+D_{l}\bar{A_{l}}z_{\tau(l)}\bar{z}_{\sigma(l)}).
\end{align*}
Since $(|A_{l}|+|D_{l}|)^{2}$ is constant, minimizing $\sum_{l\in\mathcal{L}}Q_{l}^{2}$
is equivalent to maximizing 
\begin{align}
 & \sum_{l\in\mathcal{L}}(2|A_{l}||D_{l}|+\bar{D_{l}}A_{l}\bar{z}_{\tau(l)}z_{\sigma(l)}+D_{l}\bar{A_{l}}z_{\tau(l)}\bar{z}_{\sigma(l)})\nonumber \\
 & =\sum_{l\in\mathcal{L}}(|A_{l}||D_{l}||z_{\tau(l)}|^{2}+|A_{l}||D_{l}||z_{\sigma(l)}|^{2}\nonumber \\
 & \qquad+\bar{D_{l}}A_{l}\bar{z}_{\tau(l)}z_{\sigma(l)}+D_{l}\bar{A_{l}}z_{\tau(l)}\bar{z}_{\sigma(l)})\nonumber \\
 & =z^{H}Wz\label{eq:W_quad}
\end{align}
where $z\in\C^{n}$ is the vector of variables $z_{j}$, and $W\in\C^{n\times n}$
is a Hermitian matrix whose elements are given by: \begin{subequations}\label{eq:W}
\begin{align}
W_{j,j}= & \sum_{l\in\mathcal{L}:\tau(l)=j}|A_{l}||D_{l}|+\sum_{l\in\mathcal{L}:\sigma(l)=j}|A_{l}||D_{l}|\\
W_{j,k}= & \sum_{l\in\mathcal{L}:\tau(l)=j,\sigma(l)=k}\bar{D_{l}}A_{l}+\sum_{l\in\mathcal{L}:\tau(l)=k,\sigma(l)=j}D_{l}\bar{A_{l}}\nonumber \\
 & \text{ for }j\neq k.
\end{align}
\end{subequations}
\begin{lem}
\label{lm:W_PSD} The matrix $W$ is positive semidefinite. 
\end{lem}
\begin{IEEEproof}
For every $z\in\C^{n}$, it follows from \eqref{eq:W_quad} that 
\begin{align*}
z^{H}Wz & =\sum_{l\in\mathcal{L}}(|A_{l}||D_{l}||z_{\tau(l)}|^{2}+|A_{l}||D_{l}||z_{\sigma(l)}|^{2}\\
 & +\bar{D_{l}}A_{l}\bar{z}_{\tau(l)}z_{\sigma(l)}+D_{l}\bar{A_{l}}z_{\tau(l)}\bar{z}_{\sigma(l)}).
\end{align*}
In addition, for every link $l$ it holds that 
\begin{align*}
 & \bar{D_{l}}A_{l}\bar{z}_{\tau(l)}z_{\sigma(l)}+D_{l}\bar{A_{l}}z_{\tau(l)}\bar{z}_{\sigma(l)}\\
 & =2\text{Re}(\bar{D_{l}}A_{l}\bar{z}_{\tau(l)}z_{\sigma(l)})\geq-2|A_{l}||D_{l}||z_{\tau(l)}||z_{\sigma(l)}|.
\end{align*}
Therefore, 
\begin{align*}
z^{H}Wz\geq & \sum_{l\in\mathcal{L}}(|A_{l}||D_{l}||z_{\tau(l)}|^{2}+|A_{l}||D_{l}||z_{\sigma(l)}|^{2}\\
 & \qquad-2|A_{l}||D_{l}||z_{\tau(l)}||z_{\sigma(l)}|\\
= & \sum_{l\in\mathcal{L}}|A_{l}||D_{l}|(|z_{\tau(l)}|-|z_{\sigma(l)}|)^{2}\geq0.
\end{align*}
This concludes that $W$ is positive semidefinite. 
\end{IEEEproof}
Now, one can formulate the offset optimization problem \eqref{eq:problem}
as the following QCQP: 
\begin{align}
\underset{z\in\C^{n}}{\text{maximize }} & z^{H}Wz\label{eq:QCQP}\\
\text{subject to } & |z_{j}|^{2}=1,\quad j=1,2,\dots,n.\nonumber 
\end{align}
Given a solution $\hat{z}$ to the QCQP \eqref{eq:QCQP}, one can
obtain the optimal offsets of the traffic network via the equation
\begin{align}
\theta_{s}=\frac{1}{2\pi}(\angle\hat{z}_{s}-\angle\hat{z}_{n})\label{eq:offset_angle}
\end{align}
for every intersection $s\in\mathcal{S}$.
\begin{rem}
Note that the QCQP \eqref{eq:QCQP} formulated in this paper is subtly
different from the one considered in \cite{coogan2017offset}. Specifically,
the diagonal elements of the matrix $W$ in \cite{coogan2017offset}
are all zero so the matrix is not positive semidefinite. In our formulation,
the matrix $W$ in \eqref{eq:QCQP} is positive semidefinite, which
will enable us to compute the approximation ratio of the relaxation. 
\end{rem}

\section{Approximation Algorithm\label{sec:Approximation-Algorithm}}

In the previous section, offset optimization was cast as the optimization
problem \eqref{eq:QCQP} that maximizes a convex objective function
subject to nonconvex constraints. This QCQP formulation results in
a nonconvex optimization problem. In fact, such nonconvex QCQP is
known to be NP-hard \cite{so2007approximating}. Unless P=NP, we have
to focus on finding an efficient approximation algorithm with polynomial
complexities for large-scale traffic networks.

Note that this formulation of offset optimization has a similar structure
as the QCQP formulation of the classic MAX-CUT problem in combinatorial
optimization \cite{karp1972reducibility}. Indeed, if the variable
$z$ in problem \eqref{eq:QCQP} is forced to be real, as in $z\in\R^{n},$
then the constraint $|z_{j}|=1$ implies $z_{j}\in\{+1,-1\}$, and
the maximization of a quadratic form subject to $\pm1$ variables
is exactly MAX-CUT. Consequently, we may view \eqref{eq:QCQP} as
a complex version of the MAX-CUT problem.

Based on the celebrated Goemans\textendash Williamson algorithm \cite{goemans1995improved}
for MAX-CUT, we provide below a polynomial complexity algorithm that
solves \eqref{eq:QCQP} with a performance guarantee of $\pi/4\ge0.785$
(i.e., the value of the solution is at least a factor $\pi/4$ times
the globally optimal value). In practice, the proposed algorithm might
perform even better than the provable guarantees. Our numerical results
in Section \ref{sec:numerical} find that every solution enjoys a
performance guarantee of more than $0.99$.

Following the idea of the Goemans\textendash Williamson algorithm,
one can interpret \eqref{eq:QCQP} as an optimization problem over
the one-dimensional unit sphere. This means that the problem restricts
each decision variable $z_{j}\in\C$ to be an one-dimensional unit
vector. Replacing each one-dimensional vector $z_{j}\in\C$ by an
$n$-dimensional unit vector $v_{j}\in\C^{n}$ leads to the relaxation:
\begin{align}
\underset{v_{1},\ldots,v_{n}\in\C^{n}}{\text{maximize }} & \sum_{j=1}^{n}\sum_{k=1}^{n}W_{j,k}v_{j}^{H}v_{k}\label{eq:maxcut_rel}\\
\text{subject to } & \|v_{j}\|^{2}=1,\quad j=1,\ldots,n.\nonumber 
\end{align}
This nonconvex problem can be reformulated into a convex problem by
a change of variables $X=[v_{j}^{H}v_{k}]\in\C^{n\times n}$: 
\begin{align}
\underset{X\in\C^{n\times n}}{\text{maximize }} & \tr(WX)\label{eq:maxcut_sdp}\\
\text{subject to } & X_{j,j}=1,\quad j=1,\ldots,n,\nonumber \\
 & X\succeq0.\nonumber 
\end{align}

\begin{lem}
\label{lm:relaxation} Problem \eqref{eq:maxcut_sdp} is a relaxation
of \eqref{eq:QCQP}, and therefore, its value gives an upper-bound
for the optimal value of \eqref{eq:QCQP}. 
\end{lem}
\begin{IEEEproof}
Given any feasible solution $z\in\C^{n}$ of \eqref{eq:QCQP}, let
$v_{j}=(z_{j},0,0,\ldots,0)\in\C^{n}$ for $j=1,2,\ldots,n$. Then,
$v_{j}^{H}v_{k}=\bar{z}_{j}z_{k}$ for all $j,k=1,2,\ldots,n$. Consequently,
$(v_{1},v_{2},\dots,v_{n})$ is feasible for \eqref{eq:maxcut_sdp}
and its objective value in \eqref{eq:maxcut_sdp} is the same as the
objective value of $z$ in \eqref{eq:QCQP}. 
\end{IEEEproof}
Problem \eqref{eq:maxcut_sdp} is an SDP for which an interior-point
method is able to compute an optimal solution $\hat{X}$ in polynomial
time with a given accuracy. We can recover a corresponding globally-optimal
set of vectors $\hat{v}_{1},\ldots,\hat{v}_{n}\in\C^{n}$ for \eqref{eq:maxcut_rel}
by factoring $\hat{X}=\hat{V}^{H}\hat{V}$ and taking each $\hat{v}_{j}$
to be the $j$-th column of the matrix $\hat{V}$.
\begin{rem}
The SDP \eqref{eq:maxcut_sdp} can also be generated from \eqref{eq:QCQP}
using a standard SDP relaxation procedure \cite{luo2010semidefinite}.
Specifically, by adding a rank constraint $\text{rank}(X)=1$ in \eqref{eq:maxcut_sdp},
one obtain the original QCQP \eqref{eq:QCQP} because any rank-one
matrix $X$ can be factored into $X=zz^{H}$. The relaxation \eqref{eq:maxcut_sdp}
becomes exact if its solution $\hat{X}$ has rank one. This special
situation occurs for certain types of networks~\cite{sojoudi2014exactness}
and the offsets obtained from the SDP solution achieves optimal performance
for these cases \cite{coogan2017offset}. In general, however, the
solution $\hat{X}$ of \eqref{eq:maxcut_sdp} has a rank strictly
greater than one. Nevertheless, we observe in our numerical experiments
in Section~\ref{sec:numerical} that the associated performance guarantee
(i.e. the ratio between the upper- and lower-bounds on the performance)
exceeds $99\%$ for every case. 
\end{rem}
In spirit of the Goemans\textendash Williamson idea method, one can
project an optimal set of vectors $\hat{v}_{1},\ldots,\hat{v}_{n}\in\C^{n}$
for \eqref{eq:maxcut_rel} back onto the one-dimensional unit sphere
in $\C$ by randomized rounding 
\begin{align}
 & s_{j}=r^{H}\hat{v}_{j},\qquad\hat{z}_{j}=s_{j}/|s_{j}|.\label{eq:z_round}
\end{align}
Here, $r\in\C^{n}$ is a random vector whose real and imaginary parts
are selected independently and identically from the $n$-dimensional
Gaussian distribution, as in 
\begin{equation}
r=r_{1}+ir_{2},\qquad r_{1},r_{2}\sim\mathcal{N}(0,I)\label{eq:rdef}
\end{equation}
where $\mathcal{N}(0,I)$ denotes the $n$-dimensional Gaussian distribution
with identity covariance matrix and zero mean.

This rounding method can be repeated with several choices of $r$,
and we select the candidate solution with the best objective value.
The follow result states that this randomization rounding offers a
remarkably high-quality solution. 
\begin{thm}
\label{thm:maxcut}Given the optimal solution $\hat{v}_{1},\ldots,\hat{v}_{n}\in\C^{n}$
for (\ref{eq:maxcut_rel}), define the candidate solution $\hat{z}\in C^{n}$
for \eqref{eq:QCQP} using (\ref{eq:z_round}) for each $\hat{z}_{j}\in\C$,
in which $r\in\C^{n}$ is selected as in (\ref{eq:rdef}). Then, 
\[
\sum_{j=1}^{n}\sum_{k=1}^{n}W_{j,k}\hat{v}_{j}^{H}\hat{v}_{k}\ge\text{opt}_{\text{QCQP}}\ge\ee\left[\hat{z}^{H}W\hat{z}
\right]\ge\frac{\pi}{4}\text{opt}_{\text{QCQP}},
\]
where $\text{opt}_{\text{QCQP}}$ is the globally optimal value of
\eqref{eq:QCQP} and $\ee\left[\cdot\right]$ is the expectation operator. 
\end{thm}
\begin{IEEEproof}
The first bound is true because \eqref{eq:maxcut_rel} is a relaxation
of \eqref{eq:QCQP} by Lemma \ref{lm:relaxation}, and the second
bound holds because $\hat{z}_{1},\ldots,\hat{z}_{n}\in\C$ is not
necessarily optimal for \eqref{eq:QCQP}. The third bound follows
from a result of \cite{so2007approximating}, noting that $W\succeq0$
from Lemma~\ref{lm:W_PSD}. 
\end{IEEEproof}
In summary, this section describes a $\pi/4$-approximation algorithm
for the QCQP \eqref{eq:QCQP} of the offset optimization problem that
comprises two key steps: 
\begin{enumerate}
\item Solve the SDP relaxation \eqref{eq:maxcut_sdp} and obtain the optimal
solution $\hat{X}\in\C^{n\times n}$; and 
\item Round $\hat{v}_{1},\ldots,\hat{v}_{n}\in\C^{n}$ into $\hat{z}_{1},\ldots,\hat{z}_{n}\in\C$
using the randomized procedure in \eqref{eq:z_round}. 
\end{enumerate}
Standard algorithms implement these two steps with a combined complexity
of $O(n^{4.5})$ time and $O(n^{2})$ memory, with the first step
dominating the overall complexity. These figures are polynomial, and
hence ``efficient'' in theory. In practice, however, they become
prohibitive for large-scale traffic networks with more than $1000$
intersections.

\section{Efficient Implementation for Sparse Networks\label{sec:Efficient-Implementation-for}}

When a traffic network is large but sparse in the sense that it has
a \emph{bounded treewidth} \cite{robertson1986graph}, we show in
this section that the approximation algorithm described in the previous
section can be implemented in near-linear $O(n^{1.5})$ time and linear
$O(n)$ memory.

In the following, we first describe the concept of tree decomposition
and use it to convert the original problem to a reduced-complexity
problem. Then, we further simplify the complexity to obtain a near-linear
time approximation algorithm for offset optimization.

\subsection{Tree Decomposition}



For a traffic network $G=(\mathcal{S}\cup\{\epsilon\},\mathcal{L})$,
the graph theoretical concepts of tree decomposition and treewidth
are defined as follows: 
\begin{defn}
A \emph{tree decomposition} of a graph $G$ of is a pair $(\I,T)$,
where $\I=\{\I_{1},\ldots,\I_{n}\}$ are $n$ subsets of nodes of
$G$, and $T$ is a tree with vertices $\I$, such that: 
\begin{enumerate}
\item (Node cover) For every node $s$ of $G$, there exists $\I_{j}\in\I$
such that $s\in\I_{j}$; 
\item (Edge cover) For every edge $l$ of $G$, there exists $\I_{k}\in\I$
such that $\sigma(l)\in\I_{k}$ and $\tau(l)\in\I_{k}$; and 
\item (Running intersection) If $s\in\I_{j}$ and $s\in\I_{k}$, then we
also have $s\in\I_{m}$ for every $\I_{m}$ that lies on the path
from $\I_{j}$ to $\I_{k}$ in the tree $T$. 
\end{enumerate}
\end{defn}
\begin{defn}[\cite{robertson1986graph}]
The \emph{width} of a tree decomposition $(\I,T)$ is $\omega-1$
where 
\begin{align}
\omega=\max_{j}\quad|\I_{j}|,\label{eq:omega}
\end{align}
i.e., the width is one less than the maximum number of elements in
any subset $\I_{k}\in\I$. The \emph{treewidth} of a network is the
minimum width amongst all tree decompositions. The network is said
to have a \emph{bounded treewidth} if its treewidth is $O(1)$, i.e.,
independent of the number of nodes $n$. 
\end{defn}
From the definition, the empty graph has treewidth of zero, and tree
and forest graphs have treewidths of one. Basically, the treewidth
of a graph indicates how ``tree-like'' the graph is. The treewidth
can be viewed as a sparsity criterion which determines the complexities
of many problems related to a graph. The problem of computing the
exact treewidth of a graph is known to be NP-complete~\cite{arnborg1987complexity}.
For bounded treewidth networks known \textit{a priori} to have small
$\omega\ll n$, the treewidth and the corresponding tree decomposition
can be determined in $O(2^{\omega}n)$ time~\cite{bodlaender1995approximating}.
In practice, it is much easier to compute a ``good-enough'' tree
decomposition with a small but suboptimal value of $\omega$, using
one of the heuristics originally developed for the fill-reduction
problem in numerical linear algebra. In our implementation, we use
the simple approximate minimum degree algorithm in generating a tree
decomposition \cite{amestoy1996approximate}. This approximately coincides
with the simple ``greedy algorithm'', and does not typically enjoy
strong guarantees. Regardless, the algorithm is extremely fast, generating
permutations for graphs containing millions of nodes and edges in
a matter of seconds.

Algebraically, a tree decomposition of our traffic network can also
be described by a \emph{fill-reducing permutation} matrix $P$. More
specifically, given a permutation matrix $P\in\R^{n\times n}$, we
can factor the matrix $W$ of the network into a Cholesky factor $L$
satisfying 
\begin{equation}
LL^{H}=PWP^{H},\quad L\text{ is lower-triangular},\quad L_{j,j}\ge0.\label{eq:Ldef}
\end{equation}
Let $\I_{1},\ldots,\I_{n}\subseteq\{1,\ldots,n\}$ be the column index
sets from the sparsity pattern of $L$ defined by 
\begin{equation}
\I_{j}=\{k\in\{1,\ldots,n\}:L_{k,j}\ne0\}.\label{eq:col_set}
\end{equation}
From the column index sets $\I_{1},\ldots,\I_{n}$, define a set of
parent pointers $p:\{1,\ldots,n\}\to\{1,\ldots,n\}$: 
\begin{equation}
p(j)=\begin{cases}
j & |\I_{j}|=1,\\
\min_{i}\{i>j:i\in\I_{j}\} & |\I_{j}|>1.
\end{cases}\label{eq:etree}
\end{equation}

\begin{lem}
The collection of the column index sets $\I=\{\I_{1},\ldots,\I_{n}\}$
together with the tree $T$ constructed by nodes $\I$ and edges $\{(\I_{j},\I_{p(j)}),j=1,2,\ldots,n\}$
constitute a tree decomposition for the network $G$. 
\end{lem}
\begin{IEEEproof}
According to \cite{vandenberghe2015chordal}, the pair $(\I,T)$ forms
a tree decomposition of $W$. From the definition \eqref{eq:W} of
$W$, the entry $W_{j,k}$ is zero if no link connects between the
$j$-th intersection and the $k$-th intersection. Therefore, the
sparsity pattern of the matrix $W$ is the same as the traffic network
$G$. 
\end{IEEEproof}
For networks with a bounded treewidth, we are able to find a tree
decomposition whose width is $\omega=\max_{j}|\I_{j}|=O(1)$. Since
the Cholesky factor $L$ of $W$ has at most $\omega$ nonzero elements
per column, $L$ of such networks will be a sparse matrix containing
at most $O(n)$ nonzero elements.

In the case of real-world traffic networks, the graphs are almost
\emph{planar} by construction, because the vast majority of roads
do not cross without intersecting. Planar graphs with $n$ nodes have
treewidths of at most $O(\sqrt{n})$, attained by grid graphs; a tree
decomposition within a constant factor of the optimal can be explicitly
computed using the planar separator theorem and a nested dissection
ordering. Practical traffic networks tend to have treewidths possibly
much smaller than the $O(\sqrt{n})$ figure. While local networks
may resemble grids, inter-area networks interconnecting wider regions
are more tree-like. Accordingly, their treewidth is usually bounded
by the square-root of the size of the largest grid, which is relatively
small even for networks typically thought of as ``grid-like'' such
as Manhattan and Downtown Los Angeles.

\subsection{Clique Tree Conversion and Recovery}

Using the concept of tree decomposition, this subsection describe
the clique tree conversion technique of \cite{fukuda2001exploiting}
to simplify the $\pi/4-$approximation algorithm proposed in the previous
section.

Suppose that the network has a bounded treewidth and we are given
a tree decomposition with $\omega=O(1)$ represented by a fill-reducing
permutation $P$, its associated index sets $\I_{1},\ldots,\I_{n}$,
and the parent pointers $p$. From now on, without loss of generality,
we assume that $P=I$; otherwise, we can solve the permuted problem
with $\tilde{W}=PWP^{T}$, and reverse the ordering $z=P^{T}\tilde{z}$
once a solution $\tilde{z}$ has been computed.

Given the tree decomposition, the clique tree conversion technique
reformulates (\ref{eq:maxcut_sdp}) into a reduced-complexity problem
with the variables $X_{j}\in\C^{|\I_{j}|\times|\I_{j}|},j=1,\ldots,n$:
\begin{align}
\underset{X_{1},\ldots,X_{n}}{\text{minimize }} & \sum_{j=1}^{n}\tr(W_{j}X_{j})\label{eq:ctc}\\
\text{subject to } & (X_{j})_{k,k}=1,\quad j=1,\ldots,n,\quad k\ =1,\ldots,|\I_{j}|\nonumber \\
 & R_{p(j),j}(X_{j})=R_{j,p(j)}(X_{p(j)}),\nonumber \\
 & X_{j}\succeq0,\quad j=1,\ldots,n,\nonumber 
\end{align}
where $W_{1},\ldots,W_{j}$ are matrices satisfying 
\[
\sum_{j=1}^{n}\tr(W_{j}X_{\I_{j},\I_{j}})=\tr(WX)
\]
with respect to the original $W$ matrix, over all Hermitian choices
of $X\in\C^{n\times n}$. The exact method to construct $W_{1},\ldots,W_{j}$
can be found in \cite{zhang2017sparse}. The linear operator $R_{k,j}:\C^{|\I_{j}|\times|\I_{j}|}\to\C^{|\I_{k}|\times|\I_{k}|}$
is defined to output the overlapping elements of two principal submatrices
indexed by $\I_{k}$ and $\I_{j}$, given the latter as the argument:
\begin{equation}
R_{k,j}(X_{\I_{j},\I_{j}})=X_{\I_{k}\cap\I_{j},\I_{k}\cap\I_{j}}=R_{j,k}(X_{\I_{k},\I_{k}}).\label{eq:overlap}
\end{equation}
The associated constraints $R_{p(j),j}(X_{j})=R_{j,p(j)}(X_{p(j)})$
in (\ref{eq:ctc}) are known as the \emph{overlap constraints}.

From the bounded treewidth property, this conversion reduces the number
of decision variables from $O(n^{2})$ for $X$ in \eqref{eq:maxcut_sdp}
to $O(\omega^{2}n)$ for $\{X_{j},j=1,\ldots,n\}$ in \eqref{eq:ctc}.
\begin{lem}
The solutions $\hat{X}_{1},\hat{X}_{2},\ldots,\hat{X}_{n}$ of \eqref{eq:ctc}
are related to the solution $\hat{X}$ of \eqref{eq:maxcut_sdp} by
\[
\hat{X}_{\I_{j},\I_{j}}=\hat{X}_{j},\quad j=1,\ldots,n.
\]
\end{lem}
\begin{IEEEproof}
The proof is omitted as its essentially the same as the real-valued
version in \cite{fukuda2001exploiting}. 
\end{IEEEproof}
The above relation allows us to recover the solution $\hat{X}$ of
\eqref{eq:maxcut_sdp} from solutions $\hat{X}_{1},\ldots,\hat{X}_{n}$
\eqref{eq:ctc}. Note that $\hat{X}$ is generally a dense matrix,
so simply forming the matrix would push the overall complexity up
to quadratic $O(n^{2})$ time and memory. Fortunately, the Cholesky
factorization of $\hat{X}$ is sparse due to the bounded treewidth
property. Therefore, we compute $\hat{X}$ implicitly in factorized
form a sparse factored form $\hat{X}=F^{-H}DF^{-1}$, where $D$ is
diagonal and $F$ is lower-triangular with the same sparsity pattern
as $L$ in (\ref{eq:Ldef}). This can be done by the following Algorithm~\ref{alg:recover}
adopted from \cite{vandenberghe2015chordal}. 
\begin{algorithm}[H]
\caption{\label{alg:recover}Positive semidefinite matrix completion}
\textbf{Input.} The column index sets $\I_{1},\ldots,\I_{n}$ defined
in (\ref{eq:col_set}) and the solutions $\hat{X}_{1},\ldots,\hat{X}_{n}$
to (\ref{eq:ctc}).

\textbf{Output.} The solution $\hat{X}$ to \eqref{eq:maxcut_sdp}
in the form of $\hat{X}=F^{-H}DF^{-1}$, where $D$ is a diagonal
matrix and $F$ is a lower-triangular matrix with the same sparsity
pattern as $L$.

\textbf{Algorithm.} Iterate over $j\in\{1,\ldots,n\}$ in any order.
Set $F_{j,j}=1$ and solve for the $j$-th column of $D$ and $F$
by finding any $D_{j,j}$ and $F_{\I_{j}\backslash\{j\},j}$ that
satisfy 
\[
\hat{X}_{j}\begin{bmatrix}1\\
F_{\I_{j}\backslash\{j\},j}
\end{bmatrix}=\begin{bmatrix}D_{j,j}\\
0
\end{bmatrix}.
\]
\end{algorithm}

With the factorized solution $\hat{X}=F^{-H}DF^{-1}$, we can now
efficiently implement the randomized rounding procedure described
earlier in \eqref{eq:z_round}. Specifically, from $\hat{X}=F^{-H}DF^{-1}=\hat{V}^{H}\hat{V}$
we obtain $\hat{V}=D^{1/2}F^{-1}$. Then, \eqref{eq:z_round} is equivalent
to 
\begin{equation}
F^{H}s=D^{1/2}r,\qquad\hat{z}_{j}=s_{j}/|s_{j}|\label{eq:round}
\end{equation}
where we recall that the real and imaginary parts of the random vector
$r\in\C^{n}$ are selected independently and identically from the
$n$-dimensional Gaussian distribution. Since $F$ is a lower-triangular
matrix (with the same sparsity pattern as $L$), one can compute $\hat{z}$
from \eqref{eq:round} by solving a sparse triangular system of equations
in $O(\omega n)$ time.

In summary, this subsection presents a reduced-complexity implementation
of a $\pi/4$-approximation algorithm for the QCQP \eqref{eq:QCQP}
of the offset optimization problem given a tree decomposition with
$\omega=O(1)$. The main steps are described as follows: 
\begin{enumerate}
\item Reformulate \eqref{eq:maxcut_sdp} into the reduced complexity problem
\eqref{eq:ctc}.
\item Solve \eqref{eq:ctc} to obtain solutions $\hat{X}_{1},\ldots,\hat{X}_{n}$. 
\item Recover the solution of \eqref{eq:maxcut_sdp} in the sparse factored
form $\hat{X}=F^{-H}DF^{-1}$ using Algorithm~\ref{alg:recover}. 
\item Recover a choice of $\hat{z}_{1},\ldots,\hat{z}_{n}\in\C$ via the
randomized rounding method \eqref{eq:round}. 
\end{enumerate}
We will show later that the complexity of the overall algorithm is
dominated by Step 2, i.e., the cost of solving the semidefinite program
\eqref{eq:ctc}. An interior-point method solves \eqref{eq:ctc} in
$O(\sqrt{n}$) iterations, with the cost of each iteration dominated
by the solution of a set of linear equations over $O(n)$ variables.
These equations can be fully dense despite sparsity in the original
problem, so the worst-case complexity of an interior-point solution
of \eqref{eq:ctc} is $O(n^{3.5})$ time and $O(n^{2})$ memory. Next,
we show that these complexity figures can be reduced to linear by
using dualization to exploit sparsity.

\subsection{Dualization}

Zhang and Lavaei~\cite{zhang2017sparse,zhang2017sparse2} recently
showed that the complexity of solving the real-valued version of \eqref{eq:ctc}
can be significantly improved to near-linear $O(n^{1.5})$ time and
linear $O(n)$ memory complexities via the \emph{dualization} procedure
of Löfberg~\cite{lofberg2009dualize}. We present in this subsection
a complex-valued version of the same algorithm for the traffic offset
optimization problem.

To solve \eqref{eq:ctc}, we begin by putting \eqref{eq:ctc} into
primal canonical form: 
\begin{align}
\underset{x_{1},\ldots,x_{n}}{\text{minimize }}\quad & \sum_{j=1}^{n}\bar{w}_{j}^{H}x_{j}\label{eq:vec_ctc}\\
\text{subject to }\quad & \begin{bmatrix}N_{11} & \cdots & N_{1n}\\
 & \ddots\\
N_{n1} & \cdots & N_{nn}\\
M_{1} &  & \mathbf{0}\\
 & \ddots\\
\mathbf{0} &  & M_{n}
\end{bmatrix}\begin{bmatrix}x_{1}\\
\vdots\\
x_{n}
\end{bmatrix}=\begin{bmatrix}\mathbf{0}\\
\one\\
\vdots\\
\one
\end{bmatrix},\nonumber \\
 & x_{j}\in\K_{j}\subset\C^{|\I_{j}|^{2}},\quad j=1,\ldots,n.\nonumber 
\end{align}
Each variable $x_{j}=\vect(X_{j})$ (respectively, $w_{j}=\vect(W_{j})$)
is the vectorization of $X_{j}$ (respectively, $W_{j}$) and each
$\K_{j}$ is the corresponding positive semidefinite cone. The matrices
$N_{jk}$ implement the overlap constraints in \eqref{eq:overlap}.
That is, for each $j$, the $j$-th block row $N_{j1},\ldots,N_{jn}$
implements the overlap constraint between $\I_{j}$ and its parent
$\I_{p}(j)$. Therefore, the $j$-th block row has at most two nonzero
sub-blocks: $N_{jk}=0$ except $k=j$ or $k=p(j)$. Each constraint
matrix $M_{j}$ isolates the diagonal of $X_{j}$, as in $(M_{j}x_{j})_{k}=(X_{j})_{k,k}$.

Let $N$ and $M$ denote the matrices for the constraints: 
\begin{align*}
 & N=\begin{bmatrix}N_{11} & \cdots & N_{1n}\\
 & \ddots\\
N_{n1} & \cdots & N_{nn}
\end{bmatrix}, & M=\begin{bmatrix}M_{1} &  & 0\\
 & \ddots\\
0 &  & M_{n}
\end{bmatrix}.
\end{align*}
Then, the complexity of each step of the interior-point iteration
solving \eqref{eq:vec_ctc} depends on the sparsity pattern of $\tilde{M}\tilde{M}^{H}$
where $\tilde{M}=[N^{H},M^{H}]^{H}$. Despite the nice sparsity structure
of $\tilde{M}$, the matrix $\tilde{M}\tilde{M}^{H}$ is generally
dense (see \cite{zhang2017sparse} for an example). Therefore, it
takes $O(n^{3.5})$ time and $O(n^{2})$ memory to solve \eqref{eq:vec_ctc}
using an interior-point solver.

On the other hand, the matrix $\tilde{M}^{H}\tilde{M}$ is sparse
from the block sparsity structure of $N$ and $M$. 
\begin{lem}
\label{lm:M_sparsity} The matrix $\tilde{M}^{H}\tilde{M}$ has $O(\omega^{4}n)$
nonzero elements, and it takes $O(\omega^{6}n)$ operations to compute
$\tilde{M}^{H}\tilde{M}$ from $N$ and $M$. 
\end{lem}
\begin{IEEEproof}
This is a corollary of the result of \cite{zhang2017sparse}. In particular,
$M$ is the adjacency matrix of an empty graph, so the block sparsity
structure of $\tilde{M}^{H}\tilde{M}$ is the same as the sparsity
of the adjacency matrix of the tree $T$ of the tree decomposition.
Then, $\tilde{M}^{H}\tilde{M}$ has $O(n)$ nonzero blocks, and each
of the blocks has at most $O(\omega^{4})$ nonzero elements. The computation
of $\tilde{M}^{H}\tilde{M}$ is done by adding up $O(\omega^{2}n)$
sets of blocks with $O(\omega^{4})$ elements which takes $O(\omega^{6}n)$
operations. 
\end{IEEEproof}
To exploit the sparsity structure of $\tilde{M}^{H}\tilde{M}$, Zhang
and Lavaei~\cite{zhang2017sparse,zhang2017sparse2} suggest applying
the \textit{dualization} technique of Löfberg~\cite{lofberg2009dualize}.
The main idea is to take problem \eqref{eq:vec_ctc} as it is currently
stated in primal canonical form, and rewrite the \emph{same problem}
in dual canonical form:
\begin{align}
\underset{y_{1},\ldots,y_{n}}{\text{maximize }}\quad & -\sum_{j=1}^{n}\bar{w}_{j}^{H}y_{j}\label{eq:vec_ctc_d}\\
\text{subject to }\quad & \begin{bmatrix}N_{11} & \cdots & N_{1n}\\
 & \ddots\\
N_{n1} & \cdots & N_{nn}\\
M_{1} &  & \mathbf{0}\\
 & \ddots\\
\mathbf{0} &  & M_{n}
\end{bmatrix}\begin{bmatrix}y_{1}\\
\vdots\\
y_{n}
\end{bmatrix}+s_{0}=\begin{bmatrix}0\\
\one\\
\vdots\\
\one
\end{bmatrix},\nonumber \\
 & -y_{j}+s_{j}=0,\quad j=1,\ldots,n\nonumber \\
 & s_{0}\in\{0\}^{q+1},\quad s_{j}\in\K_{j}\subset\C^{|\I_{j}|^{2}}.\nonumber 
\end{align}
Here, $q$ is the total number of equality constraints in the original
problem \eqref{eq:vec_ctc}, and $\{0\}^{q+1}$ denotes the so-called
``equality-constraint cone'', whose dual cone is a free variable
of dimension $\nu+1$.

Since \eqref{eq:vec_ctc_d} is the dual problem, with a general-purpose
primal-dual interior-point method like SeDuMi, SDPT3, and MOSEK. Each
iteration involves solving a normal equation whose block sparsity
pattern coincides with that of $\tilde{M}^{H}\tilde{M}$. Zhang and
Lavaei~\cite{zhang2017sparse2} proved that this can be done in time
linear to the number of blocks $n$, and cubic to the maximum size
of the individual blocks $\omega^{2}$. Combined, we have the following
complexity result.
\begin{thm}
\label{thm:lin_time}A general-purpose interior-point method solves
the SDP \eqref{eq:ctc} by solving its dual canonical form \eqref{eq:vec_ctc_d}
to $\epsilon$-accuracy in 
\[
O(\omega^{6.5}n^{1.5}\log\epsilon^{-1})\text{ time and }O(\omega^{4}n)\text{ memory}.
\]
\end{thm}
\begin{IEEEproof}
The proof in \cite{zhang2017sparse,zhang2017sparse2} for real-valued
SDPs can be adopted to prove this theorem. First, note that a general-purpose
interior-point method solves an order-$\theta$ linear conic program
posed in the canonical form to $\epsilon$-accuracy in $O(\sqrt{\theta}\log\epsilon^{-1})$
iterations. The cone in \eqref{eq:vec_ctc_d} has order $\theta=O(\omega n)$
from the construction of the tree decomposition, so the interior-point
method converges in $O(\omega^{0.5}n^{0.5}\log\epsilon^{-1})$ iterations. 

At each interior-point iteration, the complexity is dominated by the
solution of a linear system of equations whose block sparsity pattern
coincides with $\tilde{M}^{H}\tilde{M}$. From Lemma \ref{lm:M_sparsity},
forming $\tilde{M}^{H}\tilde{M}$ requires $O(\omega^{6}n)$ time
and $O(\omega^{4}n)$ memory. Applying \cite[Lemma 5]{zhang2017sparse2},
the cost of solving the linear equations is also $O(\omega^{6}n)$
time and $O(\omega^{4}n)$ memory. Combined, the memory complexity
is $O(\omega^{4}n)$, and the time complexity is $O(\omega^{6}n)$
per-iteration multiplied by $O(\omega^{0.5}n^{0.5}\log\epsilon^{-1})$
iterations. 
\end{IEEEproof}

\subsection{Overall Algorithm}

This section presents a reduced-complexity implementation of a $\pi/4$-approximation
algorithm for the QCQP \eqref{eq:QCQP} of the offset optimization
problem. The full algorithm is described as follows: 
\begin{enumerate}
\item Compute a tree decomposition for the traffic network $G$ and its
fill-reducing permutation $P$ using the minimum degree algorithm.
\item Permute $W$ as $W\gets PWP^{H}$, compute the Cholesky factor $L$
as in \eqref{eq:Ldef}, and determine the index sets $\I_{1},\ldots,\I_{n}$
and the parent pointers $p$, as in \eqref{eq:col_set} and \eqref{eq:etree}. 
\item Use the clique tree conversion technique to reformulate \eqref{eq:maxcut_sdp}
into \eqref{eq:ctc}. 
\item Convert \eqref{eq:ctc} to the dualized problem \eqref{eq:vec_ctc_d}. 
\item Solve \eqref{eq:vec_ctc_d} as a dual canonical problem using a general-purpose
primal-dual interior-point method to obtain solutions $\hat{X}_{1},\ldots,\hat{X}_{n}$
of \eqref{eq:ctc}. 
\item Recover the solution of \eqref{eq:maxcut_sdp} in the sparse factored
form $\hat{X}=F^{-H}DF^{-1}$ using Algorithm~\ref{alg:recover}. 
\item Recover a choice of $\hat{z}_{1},\ldots,\hat{z}_{n}\in\C$ via the
randomized rounding method \eqref{eq:round}. This randomization step
can be run several times to obtain a solution with the best objective
value.
\item Reverse the fill-reducing permutation $\hat{z}\gets P^{H}\hat{z}$. 
\end{enumerate}
\begin{cor}
\label{cor:complexity} The proposed algorithm generates a choice
of $\hat{z}_{1},\ldots,\hat{z}_{n}\in\C$ that satisfy the bounds
in Theorem~\ref{thm:maxcut} and can be computed with the same time
and memory complexity as described in Theorem~\ref{thm:lin_time}. 
\end{cor}
\begin{IEEEproof}
The minimum degree algorithm in Step 1 takes $O(\omega n)$ time and
memory. Step 2 is dominated by the Cholesky factorization step, for
$O(\omega^{3}n)$ time and $O(\omega^{2}n)$ memory. Steps 3 and 4
are algebraic manipulations, requiring $O(\omega^{2}n)$ time and
memory. Step 5 uses $O(\omega^{6.5}n^{1.5}\log\epsilon^{-1})$ time
and $O(\omega^{4}n)$ memory according to Theorem~\ref{thm:lin_time}.
Algorithm \ref{alg:recover} in Step 6 is dominated by solving $n$
linear systems of up to size $\omega^{2}$ for $O(\omega^{3}n)$ time
and $O(\omega^{2}n)$ memory. The rounding method of \eqref{eq:round}
in Step 7 can be performed by back-substitution in $O(\omega n)$
time and memory. Finally, Step 8 takes $O(n)$ time and memory to
obtain an approximate solution with the guarantees in Theorem~\ref{thm:maxcut}. 
\end{IEEEproof}
\begin{rem}
The offset optimization problem \eqref{eq:QCQP} is formulated as
a complex-valued QCQP. This complex-valued QCQP has an equivalent
real-valued formulation. Specifically, consider $z=x-iy$ where $x,y\in\R^{n}$
are the real and imaginary parts of $z$. Then, \eqref{eq:QCQP} is
equivalent to 
\begin{align}
\underset{x,y\in\R^{n}}{\text{maximize }} & [x^{H}y^{H}]\begin{bmatrix}\re{W} & \img{W}\\
-\img{W} & \re{W}
\end{bmatrix}\begin{bmatrix}x\\
y
\end{bmatrix}\label{eq:real_QCQP}\\
\text{subject to } & x_{j}^{2}+y_{j}^{2}=1,\quad j=1,2,\dots,n.\nonumber 
\end{align}
One can then follow a similar procedure to solve this transformed
real-valued problem as in our conference version \cite{ouyang2018cdc}.
However, transforming a complex QCQP into its real-valued counterpart
also doubles its treewidth. In practice, the resulting algorithm is
about a constant factor of $10$ times slower than the one proposed
in this paper. See \cite{cedric_complex} for such speed-up in optimization
solvers using complex numbers instead of real numbers. 
\end{rem}

\section{Numerical Experiments}

\label{sec:numerical}

In the previous sections, we proved that our algorithm solves offset
optimization to a global optimality ratio of $\pi/4\ge0.785$ in near-linear
$O(n^{1.5})$ time. In this section, we benchmark these guarantees
on two datasets:
\begin{enumerate}
\item Real-world dataset for the Manhattan network, based on real network
topology, flow rates and turning ratios.
\item Synthetic dataset for the Berkeley, Manhattan, and Los Angeles networks,
with real network topologies but synthetic flow rates and turning
ratios.
\end{enumerate}
In our numerical results described below, the empirical time complexity
of the algorithm is linear $O(n)$, and the computed solutions have
global optimality ratios exceeding $0.99$.

\subsection{Real-world Manhattan dataset}

We demonstrate our algorithm in a real-world setting, by solving offset
optimization on a traffic model of Manhattan from Osorio et al.~\cite{osorio2015scalable}
based on real data. Our network graph contains 189 nodes and 472 edges,
and covers the area between 7th and 12th Avenues, and 30th and 50th
Streets. Detailed traffic simulations were performed to result in
five sets of flow rates and turn ratios. In each case, green splits
were assigned in order to make north-south links completely out of
phase with east-west links.

We implement our algorithm in MATLAB and perform our experiments on
a 3.3 GHz 4-core Intel Xeon E3-1230 v3 CPU with 16 GB of RAM. For
each set of flow rates and turn ratios, we solve the convex relaxation
to obtain a lower-bound (``lower''), and perform randomized rounding
200 times to obtain a suboptimal solution and an upper-bound (``upper'').
As shown in Table~\ref{table0}, all of the five optimization problems
completed within ten seconds, to result in global optimality ratios
of $\ge0.996$. We emphasize that global optimality must be interpretted
within the context of the formulation. In particular, they assume
that traffic flow can be adequately approximated as being sinusoidal. 

\begin{table}[t]
\caption{ Offset optimization on the real-world Manhattan dataset}
\label{table0} \centering
\resizebox{\columnwidth}{!}{ %
\begin{tabular}{|c|c|c|c|c|c|c|}
\hline 
dataset & 1 & 2 & 3 & 4 & 5 & mean\tabularnewline
\hline 
upper & 93591 & 104339 & 96160 & 107935 & 98639 & 100133\tabularnewline
\hline 
lower & 93544 & 104267 & 96159 & 107740 & 98247 & 99991\tabularnewline
\hline 
ratio & 0.9995 & 0.9993 & 1.0000 & 0.9982 & 0.9960 & 0.9986\tabularnewline
\hline 
sec & 4.46 & 3.60 & 4.91 & 3.53 & 3.89 & 4.08\tabularnewline
\hline 
\end{tabular}}
\end{table}

\subsection{Synthetic OpenStreetMap dataset}

To benchmark the scalability of our algorithm over a range of network
sizes, we generate synthetic test cases using real-world network topologies
collected from the OpenStreet Map data \cite{OSM}. For each test
case, we consider a rectangular area of the real-world map. From each
area, we construct a traffic network by assuming that all intersections
in the area are signalized. Entry links are added for roads/ways entering
the target rectangular area, and a non-entry link is added from one
intersection to another one if there is a road/way between the two
intersections following the corresponding direction. We assume that
vehicles travel at a constant speed, so the travel time $\lambda_{l}$
of each link is assigned to be proportional to the length of the link
in the real-world map. The turn ratios $\beta_{lk}$'s are set to
be such that, when vehicles entering an intersection form a link,
the traffic traveling straight is twice the traffic making each turn
direction (left or right). The average flow $f_{l}$'s of all entry
links are assumed to be the same constant, and the flows of non-entry
links are calculated from the turn ratios by solving $f_{l}=\sum_{k\in\mathcal{L}}\beta_{kl}f_{k}$
for all $l\in\mathcal{L}$.

Since the focus is on the offsets, other signal control parameters
are set to be fixed. The cycle lengths of all intersections are the
same constant as described in the network model. For each network,
the splits and phase sequences are described by the green split parameters
$\gamma_{l}$. In the numerical experiments, we do not optimize the
green splits $\gamma_{l}$ and set them based on the orientations
of the links for convenience. In particular, at each intersection,
the green split of a link is the angle between the corresponding road/way
and the longitude line of the intersection on the real-world map.

\begin{figure}
\begin{centering}
\subfloat[]{\begin{centering}
\includegraphics[width=1\columnwidth]{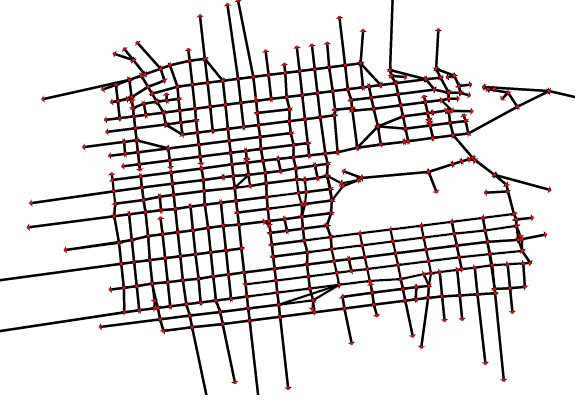} 
\par\end{centering}
}
\par\end{centering}
\begin{centering}
\subfloat[]{\begin{centering}
\includegraphics[width=1\columnwidth]{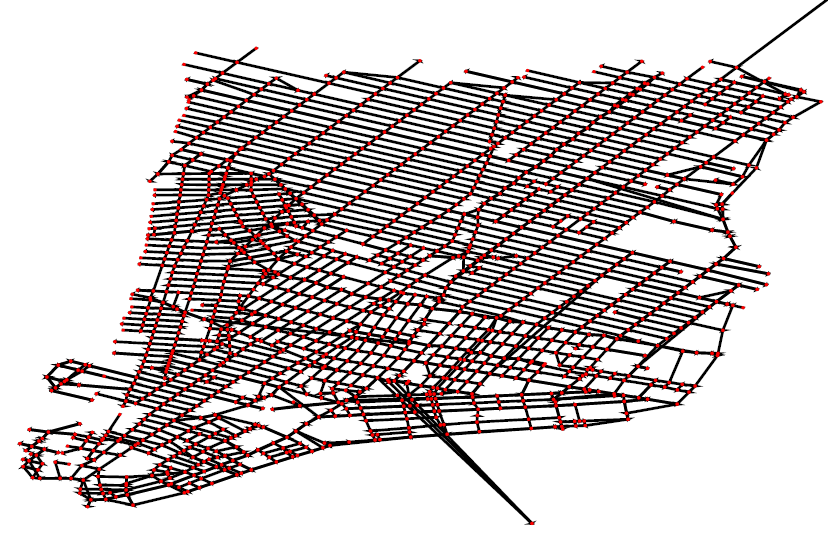} 
\par\end{centering}
}
\par\end{centering}
\begin{centering}
\subfloat[]{\begin{centering}
\includegraphics[width=1\columnwidth]{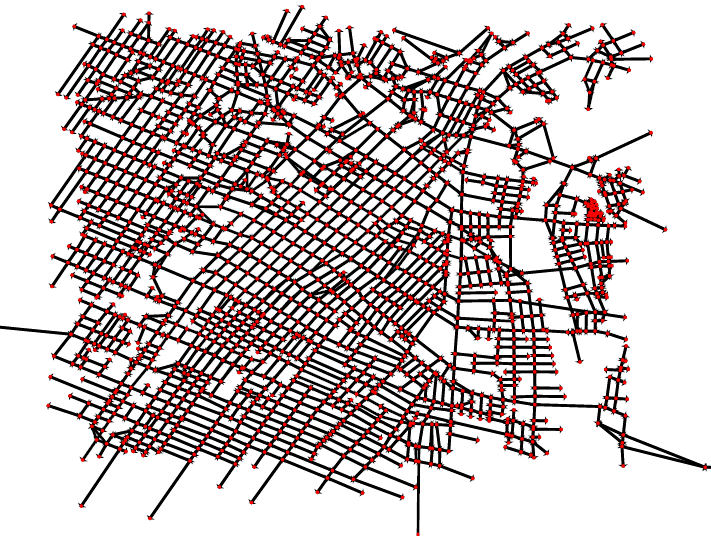} 
\par\end{centering}
}
\par\end{centering}
\caption{The three road networks under study: (a) Berkeley; (b) Manhattan;
(c) Downtown Los Angeles}
\label{fig:networks} 
\end{figure}

The first set of the networks is generated using the map of the Berkeley
area as shown in Fig.~\ref{fig:networks}a. The Berkeley-1 network
has $405$ intersections and $1122$ links connecting the intersection,
while Berkeley-4 has $12176$ intersections and $33725$ links that
includes the network of Berkeley, Oakland, and their surrounding areas.
The second set of networks is generated from the map of the Manhattan
area as in Fig.~\ref{fig:networks}b, and the third set of networks
is based on the Downtown Los Angeles area as in Fig.~\ref{fig:networks}c.

The network parameters and numerical results are presented in Table
\ref{table1}. The number of intersections $n$ ranges from $405$
to $12176$ among the networks in our experiments. In every case,
the tree decomposition parameter $\omega$ is bounded by $50$. The
lower bound (``lower'') on the squared queue length is the optimal
value of the optimization problem \eqref{eq:ctc} obtained from Step
5 of the algorithm. The optimal value of \eqref{eq:ctc} serves as
a bound according to Theorem \ref{thm:maxcut}. For each network,
the upper bound (``upper'') is the result from the best solution
$\hat{z}$ in $200$ runs of the randomized rounding method in Step
7 of the algorithm. The algorithm is implemented in MATLAB, and the
numerical experiments are performed on an HP SE1102 server with 2
quad-core 2.5GHz Xeon and 24 GB memory.

\begin{table}[t]
\caption{Offset optimization on the synthetic Berkeley (``Berk''), Manhattan
(``NYC''), and Los Angeles (``LA'') datasets}
\label{table1} \centering
\resizebox{\columnwidth}{!}{ %
\begin{tabular}{|c|c|c|c|c|c|c|c|}
\hline 
Cases  & $|\mathcal{S}|=n$  & $|\mathcal{L}|$  & $\omega$  & lower & upper & ratio & sec\tabularnewline
\hline 
Berk-1  & 405  & 1122  & 14  & 79209  & 79498  & 0.9964 & 6 \tabularnewline
\hline 
Berk-2  & 2036  & 5789  & 36  & 477449  & 479725  &  0.9953 & 253 \tabularnewline
\hline 
Berk-3  & 6979  & 19222  & 41  & 1588518  & 1597089  & 0.9946 & 1253\tabularnewline
\hline 
Berk-4  & 12176  & 33725  & 42  & 2795240  & 2810684  & 0.9945 & 2657 \tabularnewline
\hline 
NYC-1  & 1430  & 2748  & 37  & 301366  & 303057  & 0.9944 & 234\tabularnewline
\hline 
NYC-2  & 2016  & 3854  & 31  & 417186  & 419692  & 0.9940 & 232\tabularnewline
\hline 
NYC-3  & 3923  & 7841  & 37  & 780878  & 787526  & 0.9916 & 655\tabularnewline
\hline 
NYC-4  & 9968  & 20945  & 39  & 2022529  & 2039907  & 0.9915 & 2565\tabularnewline
\hline 
LA-1  & 733  & 2180  & 22  & 182811  & 183403  & 0.9968 & 28\tabularnewline
\hline 
LA-2  & 1838  & 5170  & 36  & 458209  & 460708  & 0.9946 & 171 \tabularnewline
\hline 
LA-3  & 3062  & 8838  & 43  & 747805  & 752536  & 0.9937 & 707 \tabularnewline
\hline 
LA-4  & 4239  & 12773  & 50  & 1139072  & 1146237  & 0.9937 & 2207\tabularnewline
\hline 
\end{tabular}}
\end{table}

\begin{figure}[t]
\begin{centering}
\includegraphics[width=1\columnwidth]{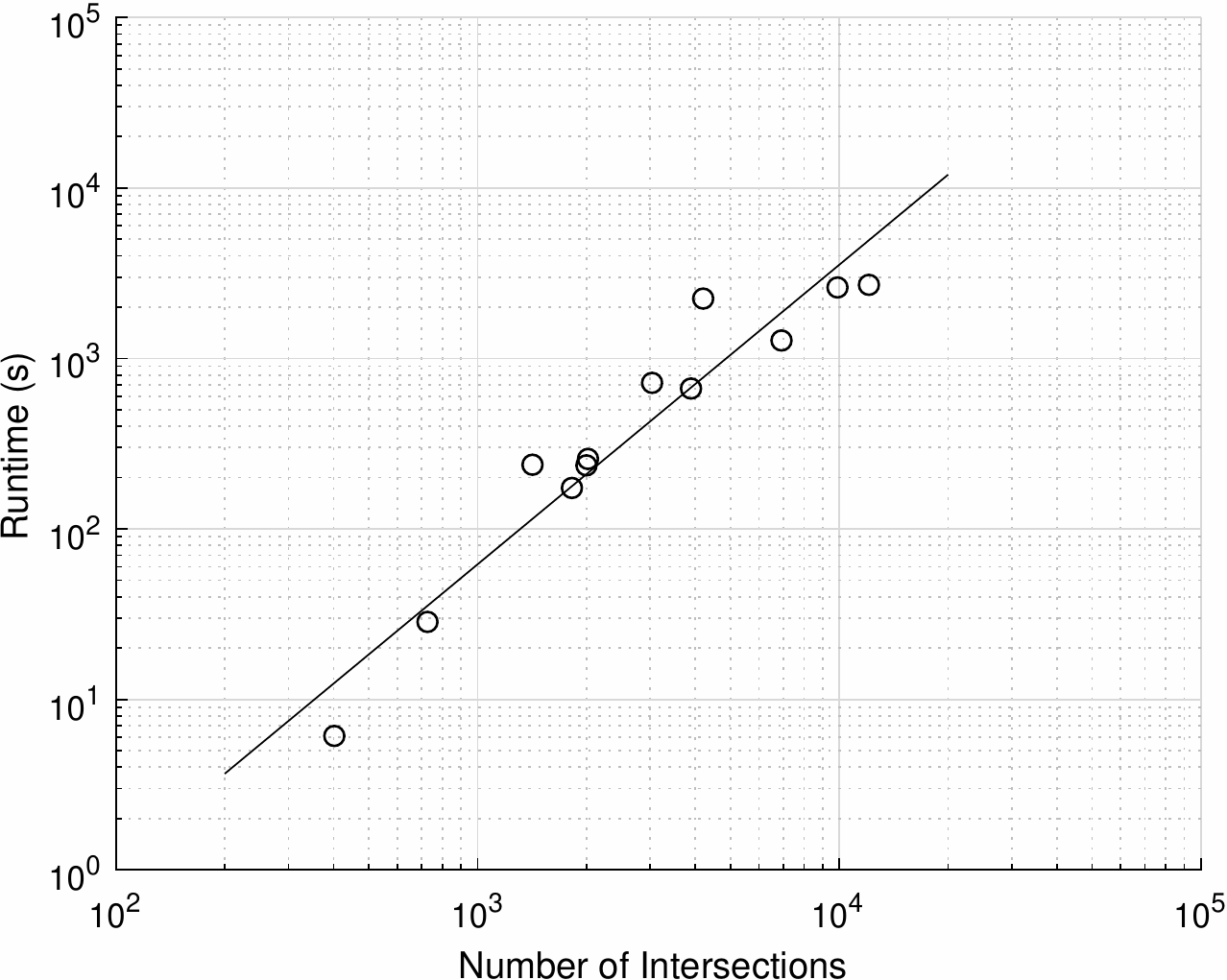} 
\par\end{centering}
\caption{Runtime against number of intersections. The regression line plots
a \emph{linear} empirical time complexity.}
\label{fig:runtime} 
\end{figure}

As observed in Table \ref{table1}, the performance of the algorithm
is much better than the theoretical (worst-case) $\pi/4$ guarantee
in Theorem \ref{thm:maxcut}. In fact, the gap between the upper and
lower bounds on the queue lengths is less than $1\%$ for all cases
($99\%$ optimal guarantee). Therefore, despite being an approximation
algorithm, the proposed algorithm is able to provide almost globally
optimal solutions for the offset optimization problem generated from
real-world traffic networks.

In terms of runtime, the algorithm can solve the SDP relaxations and
compute near-optimal offsets for networks with up to twelve thousand
intersections within an hour. This allows the potential to re-compute
offsets every hour based on real-time traffic conditions. Furthermore,
Fig. \ref{fig:runtime} shows that the runtime scales almost linearly
with respect to the number of intersections in the network. This agrees
with the claim of Corollary~\ref{cor:complexity} and it demonstrates
the ability of our algorithm in solving large-scale traffic offset
optimization. 

\section{Conclusion}

We describe an algorithm that solves the Coogan et al.~\cite{coogan2017offset}
formulation of traffic signal offset optimization to near-global optimality
in near-linear time. The algorithm performs a randomized rounding
of an SDP relaxation to yield a suboptimal solution with global optimality
bound $\pi/4\ge0.785$. Assuming that the traffic network has a \textquotedblleft tree-like\textquotedblright{}
topology, we prove that the algorithm has $O(n^{1.5})$ time complexity
and $O(n)$ memory complexity with respect to the number of intersections
$n$. Numerical experiments verify the underlying complexity result,
and the algorithm is able to obtain near-globally optimal solutions
for networks with up to twelve thousand intersections within an hour.
An important future work is to benchmark these results against realistic
large-scale micro-simulations, in order to validate the assumptions
underlying the traffic model.

\section*{Acknowledgments}

The authors are grateful to Carolina Osorio and Timothy Tay at MIT
for access to their real-world traffic models of Manhattan, and for
their invaluable help in generating meaningful datasets. 




\ifCLASSOPTIONcaptionsoff \newpage\fi




\bibliographystyle{IEEEtran}
\bibliography{../IEEEabrv,references}

\end{document}